\documentclass[twoside,12pt]{article}
\usepackage{latexsym,amssymb,theorem, epsfig, color}
\usepackage{amsfonts}
\usepackage{amsmath}
\DeclareMathSymbol{\C}{\mathalpha}{AMSb}{"43}
\DeclareMathSymbol{\R}{\mathalpha}{AMSb}{"52}
\DeclareMathSymbol{\Z}{\mathalpha}{AMSb}{"5A}

\newtheorem{theorem}{Theorem}[section]
\newtheorem{lemma}[theorem]{Lemma}
\newtheorem{definition}[theorem]{Definition}

\newtheorem{remark}[theorem]{Remark}
\newtheorem{observation}[theorem]{Observation}

\newcommand{\qed}{\nobreak \ifvmode \relax \else
      \ifdim\lastskip<1.5em \hskip-\lastskip
      \hskip1.5em plus0em minus0.5em \fi \nobreak
      \vrule height0.75em width0.5em depth0.25em\fi}

\title{Total wave based fast direct solver for VSP}
\author{Yu Chen \\
Courant Institute of Mathematical Sciences \\
New York University}
\date{Aug 18, 2010}

\begin{document}
\maketitle
\tableofcontents

\vspace{0.2in}

\noindent
{\bf Summary of the report.}
We present a fast direct solver for the volume scattering problem
of the Helmholtz equation. The algorithm is faster than existing 
methods. Moreover, discretization for our method is much simpler 
and more accurate than that for finite difference, finite elements,
and integral equations.

Jacques Hadamard's work on ill-posedness put us in a box of
solving well-posed problems of preferably small condition number. 
In reformulating elliptic problems, such as the scattering problem
for the Helmholtz or Maxwell, by integral equations (IEs) 
the price we have to pay is the complexity in solution representations
and their discretizations, and cost of computation (20 or more, 
instead of 2 to 4, points per wavelength) for a fast direct solver
(FDS). Well-posedness lead us, for example, to reformulating the 
first order Maxwell equations as second order elliptic PDEs, with 
symmetries breached and balances of physical quantities disturbed. 
Additional difficulties will have to be met when the second order 
equations are formulated by IEs, and as the latter are discretized 
and solved numerically. 

There is another all-embracing box next to that of Hadamard's. 
The more capable we become of solving problems, the more we seek
challenging and interesting problems, such as a scattering problem. 
In the present report we propose a different approach to solution 
of elliptic problems in a compactly supported domain $D$ 
with inhomogeneous medium, referred to here as the volume 
scattering problem (VSP). Our method never solves a well-posed or
ill-posed problem. Instead it first solves a very simple but 
not posed problem: PDEs without boundary conditions.

We will construct the total wave solution space (TWSS) for the 
given PDEs, or more precisely, we will construct the null space 
of the homogenous PDEs (with zero RHS and of variable coefficients),
subject to no boundary or any other conditions. Thus the TWSS consists 
of {\it general} solutions of the PDEs. Constructing the null space 
seems a non-scattering problem but it may be the easiest way to 
account for global communictions of the linear system, or global 
multiple scatterings of a scattering problem. It is only after the 
TWSS is constructed in $D$ that the data, or  
the incident wave, is incorporated to obtain a {\it specific} 
solution - the specific total wave corresponding to the specified 
incident wave. 

Our method will give rise to a fast direct solver for the elliptic
problem in $D$. It does not support an iterative solver. The TWSS 
will be constructed recursively on a (quadtree for example, in 2-D) 
hierarchy of domain decomposition, with the TWSS first constructed 
in each bottom level subdomains. Merging the total waves in the 
subdomains to those of their parents will end up with the TWSS 
constructed efficiently for the entire scatterer $D$.

The method will be presented in the scattering context and language, 
but the principles extend to a general boundary value problem for 
elliptic PDEs, in particular the Maxwell equations. There are 8 
first order equations for the 6 unknown $E$ and $B$ for linear 
material. We can construct the null space of the linear equations 
in $D$ first, and deal with the boundary conditions latter. Another 
major benefit of the total wave approach is that it greatly simplifies
discretization. 

In this report, a total wave refers to a nontrivial solution of the 
homogenous PDEs with variable coefficients in domain $D$.

\section{Introduction}  \label{sec-1}

The subject of this report is solution of scattering problem for 
the Helmholtz equation. There are two standard
types of scattering problems. One is for the inhomogeneous medium
inside a domain and is referred to here as the volume scattering 
problem (VSP). The other is for an impenetrable scatterer (such 
as a perfect conductor for Maxwell), or a penetrable scatterer 
with constant coefficients inside, and is referred to as surface 
scattering problem (SSP). When formulated as integral equations 
(IEs), VSP is related to the Lippmann-Schwinger IE, whereas SSP 
is related to boundary IEs. 

These problems can also be solved by finite element or boundary 
element methods.

There are fast direct solvers (FDSs) and iterative solver. A FDS is 
desirable when the condition number of the problem is not small, which 
occurs, for example, when the problem is near resonance. They are also 
more efficient for a scattering calculation with multiple incident
waves typically required for an inverse scattering problem. 

For IEs or finite element methods, discretization has never been 
made robust, or even easy. For a second kind IE, with all its 
underlying benefits in conditioning and reduction of dimensionalities, 
the discretization problem is even
more evident. The use of layer potentials makes their discretization
extremely unwieldy and difficult. For example, a quadrature to
integrate the Coulomb potential $1/r$ near $y$, with $r=|x-y|$ and 
with $y$ a fixed point on a patch of smooth surface, against smooth 
functions such as polynomials, is not easy to design due to the 
strong and non-trivial influence of the curvatures on $1/r$ near 
the source $y$.  Warping a smooth function we get a smooth function. 
Warping a singular function, we'd better prepared to reap the 
whirlwinds. So far we have not looked at the frightening situations
when these singular kernels meet corners and edges on surfaces or 
inside inhomogeneous media. 

The desirable analytic properties and attractiveness of IE 
formulation for the scattering problems, or for any elliptic 
PDEs, or for Maxwell equations, become difficult to exploit the
moment they meet discretization. This is because after discretization
we have to deal with individual and standing alone poles and dipoles. 
A pole will not interact well with other nearby poles or singularities 
unless they are on the same patch of smooth surface. Unfortunately, 
for many interesting applications, different parts of the surface, 
or multiple inclusions, may get very close.  On the fly, designing 
a quadrature will be more difficult than resampling.

We propose a solution method for the scattering problem, also for
other elliptic PDEs and the Maxwell equations, which does not require 
discretization of IEs or PDEs. It requires what we call ``sampling'' 
of the original PDEs (as in a collocation method) and their
solutions. We don't call it discretization in the sense that the
collocation method may be regarded as sampling, as opposed to 
discretizing the differential equation. There is a distinction 
between sampling and discretization. Sampling requires not much
of brain, whereas discretization requires too much of it.

Our method does not reformulate the original PDEs as IEs. It works 
with the PDEs directly. It barely solves those PDEs. Certainly it 
never solves a well or ill-posed problem for the PDEs. It solves a
not posed problem for the PDEs. It solves the PDEs without boundary 
or other conditions. For scattering problems or other elliptic
PDEs formulated as a boundary value problem, the boundary values will
be processed only after the {\it general} solution space for the 
PDEs without boundary conditions are constructed. Constructing the
general solution space is easier than solving boundary value problems.

The method is a fast direct solver; it cannot be related to and does 
not support any iterative solver. For a SSP, with the surface not 
very concave or convoluted, our FDS is faster than existing FDSs. 
For VSPs, our FDS offers the same asymptotic complexity, with 
a big reduction on the constant if the inhomogeneous medium occupies
a convex domain $D$, such as a square or triangle.

Organization of the report: \S \ref{sec-2} is a short and
informal description of the method. \S \ref{sec-3} is a full and
more formal description of the method. \S \ref{sec-4} provides 
formulations of the volume scattering problem and Green's third 
identity used as a projector on the boundary of a scatterer. 
\S \ref{sec-5} contains background information on layer potential
representation for the interior and exterior projectors.

\section{Our method - Informal description}  \label{sec-2}

We will present the method in the context of VSP for the Helmholtz; 
see \S \ref{sec-4} for more details on VSP. For a given precision 
$\epsilon>0$ and incident wave $u_0$, the algorithm finds the unique 
solution to the scattering problem in a compactly supported 
inhomogeneous medium (a variable index of refraction $n(x)$) inside
domain $D$. If $v$ denotes the scattered wave, then $u=u_0+v$ is the 
total wave, which satisfies the homogeneous Helmholtz with variable 
coefficients
\begin{equation}
  \Delta u + k^2 n^2(x)u = 0, \quad x \in D
\label{1.1a}
\end{equation}

\subsection{Informal description}  \label{sec-2.1}

Our method consists of three parts. 

{\bf Part I.} For the prescribed precision $\epsilon>0$, a complete 
set of solutions $\{ u_j, j=1:N \}$ to the homogeneous Helmholtz 
(\ref{1.1a}) are efficiently constructed in $D$ and made available 
on the boundary $\partial D$. The boundary lies in the free space. 
See \S \ref{sec-3} for further details on the size $N$ and how the $N$ 
solutions are constructed efficiently. 

{\bf Part II.} Each solution $u_j$ thus constructed, being a total 
wave, can and will be split into two parts on $\partial D$, the 
incident and scattered waves $u_{0j}$ and $v_j$. This can be 
accomplished with the third Green's identity used as a projector on 
$L^2(\partial D)$. Since the set of total waves $\{ u_j \}$ is complete,
any total wave, namely a solution of the homogeneous Helmholtz 
(\ref{1.1a}), can be represented by linear combination of 
$\{ u_j, j=1:N \}$ to precision $\epsilon$. Likewise, any
incident wave can be represented by linear combination of 
$\{ u_{0j}, j=1:N \}$ to precision $\epsilon$. In particular,
our prescribed incident wave $u_0$ can be expressed in terms of 
$\{ u_{0j}, j=1:N \}$
\begin{equation}
  u_0(x) = \sum_{j=1}^N c_j u_{0j}(x), \quad x \in \partial D
\label{1.2h}
\end{equation}
Solving (\ref{1.2h}) for the coefficients $c_j$, we obtain the 
scattered wave
\begin{equation}
  v(x) = \sum_{j=1}^N c_j v_j(x), \quad x \in \partial D
\label{1.2i}
\end{equation}
corresponding to $u_0$. At this point, we have obtained the scattered 
wave on $\partial D$ and consequently also outside $D$. For many 
applications, such as inverse scattering by repeatedly solving forward 
problems, the scattered wave outside the medium is all that we want.

{\bf Part III.} Now suppose we also want the scattered wave $v$
inside $D$. It will be obtained efficiently by a downward 
recursive procedure along a hierarchical structure, such as a 
quadtree for 2-D domain $D$, which was also used in Part I to 
efficiently construct the total waves $\{ u_j, j=1:N \}$ in the 
first place.

\subsection{Discussions}  \label{sec-2.2}

The hierarchical mergings for Part I, and splittings for Part III 
are universally employed in a typical FDS, although they may not 
always be presented in a familiar language or structure. It will
referred to as domain decomposition. 
\begin{definition}
\label{def-subdm}
Throughout the report, $D$ is regarded as singly connected. A
subdomain is always the result of partitioning $D$ artificially 
for the domain decomposition. 
\end{definition}

The method can also be adopted for SSP and scattering problems
in layered media, so that no layered Green's function is required 
which is necessary for the IE formulation. See \S \ref{sec-3.3} 
for more details on the extensions to SSP.

The method solves no scattering problem in order to construct the
TWSS. It never solve a well or ill-posed problem for PDEs. In fact 
one can largely roam in the null space of linear operators without 
encounter ill-posedness; see \S \ref{sec-2} for further details. 
The modern concept of first or second kind IEs seems oxymoron to 
the primitive kernel hunter-gatherer, and null spaces of IEs or
geometric resonances are his trophies to hang up on the walls, 
not to become his stumbling blocks. 

A total wave in a subdomain, such as a square, subject to no condition
on the boundary, is only aware of the medium inside the domain. It has 
no knowledge of what is outside, in particular whether the medium is
discontinuous over the boundary, and thus it is unaware of the corners 
of the domain, unless explicitly informed; see (\ref{3.31}). In 
contrast, a scattered wave is a solution of an inhomogeneous
Helmholtz, subject to outgoing radiation conditions in the free
space, in which the subdomain must be embedded to set up the scattering
problem for the subdomain, on whose solution a scattered wave based
FDS relies. Thus the scattered wave sees the manmade discontinuity across 
the boundary. It is aware of the corners and requires more points
there to be represented.

\section{Detailed description of the Algorithm}  \label{sec-3}

We first present in \S \ref{sec-3.1} the basic components required by
our fast direct solver. We then describe our TWSS based FDS in \S 
\ref{sec-3.2}. The algorithm is similar to those of 
\cite{chen-FDS1}-\cite{h1-mat1} in data structure and complexity; 
Our method differs from theirs in what solution space to construct and 
how to merge solution spaces of subdomains; see Remark \ref{rmk-diff}.

\subsection{The Basic Components and Parameters}  \label{sec-3.1}
For simplicity, we assume that the scatterer $q$ is a smooth function, 
which vanishes smoothly outside a square domain $D$. In implementation, 
this requirement can be relaxed to include piecewise smooth scatterers 
with jumps in a fairly arbitrary, bounded domain $D$ in two dimensions. 

A typical fast direct solver for the Helmholtz equation with
large wave number $k$ relies on domain decomposition of some sort
\cite{chen-FDS1}-\cite{ming-gu}. For the Lippmann-Schwinger equation 
(\ref{3.17}), the square domain $D$ is partitioned hierarchically
into the balanced quadtree; again for simplicity we will not 
discuss adaptive partitioning until \S \ref{sec-3.3}.

The size of a scattering problem is measured by the number
of wavelengths in each linear dimension, and the number of points
required to discretize $D$ is therefore proportional to $k^2$. Let
$N=O(k^2)$ be the number of unknowns in the resulting linear system
of equations to be solved. It is well known, \cite{chen-FDS1},
\cite{Chew}, that a fast direct solver requires $O(N^{1.5}) = O(k^3)$
flops to construct general solution space and additional $O(N\log N)$ 
flops to obtain the specific scattered wave for the prescribed 
incident wave. Our approach has the same complexity, with a 
considerable reduction on the constant. 

Domain decomposition can be carried out, in principle, {\it before} or 
{\it after} the discretization. Our approach will be able to cleanly 
separate the two issues in a straightforward way, and partition $D$
before discretization. In contrast, existing implementations go the 
other way around, and have to deal with ``subdomains'' of mesh points 
of the discretization. Linear algebra tricks combined with untidy local 
approximation steps are employed to copy with these ``subdomains'' and
the communications among them. The artificial cuts and corners must 
still manifest themselves in the ``subdomains of points''. For the
Lippmann-Schwinger equation, the points are monopoles or dipoles. 
Strategies have been designed and strifes directed to these issues, 
to try to mend or heal the cuts and wounds. 
\begin{observation}
\label{obs-merg1}
For a prescribed precision and in a subdomain, typically a square or 
a rectangle obtained by merging two squares, the number of distinct solutions, 
whether the incident waves, scattered waves, or total waves, in the 
subdomain is proportional to the arclength of the subdomain, as 
measured by the number of wavelengths. This is because these solutions
can be determined uniquely by their Dirichlet, or Neumann, or D\&N 
data on the boundary of the subdomain.
\end{observation}
\begin{definition}
\label{def-solusp}
 For a prescribed precision the finite number of distinct solutions 
 form the solution space for the subdomain in question.
\end{definition}
\begin{remark}
\label{rmk-arclen}
For simplicity, we will say that the dimension of solution space for 
a square subdomain, of edge length $L$ and boundary arclength $4L$, 
is $4L$, instead of proportional to $4L$. Therefore, the dimension 
of solution space for a rectangle by merging two squares is $6L$.
\end{remark}
A necessary step of a fast direct solver is to construct the solution
space for the entire $D$ efficiently in a hierarchical order. Once
solution space is formed for each subdomain on the same level of the
quadtree, two neighbor square subdomains on the level are merged together 
to form the solution space for the union of two squares, and so on. 
\begin{remark}
\label{rmk-diff}
The main difference between existing methods and our approach is that (i)
We merger the total waves whereas they merge the scattered waves (ii)
They organize the scattered waves by the incident waves whereas we don't 
organize the total waves by the incident waves or any other waves (iii)
We merge the total waves by simple continuity conditions whereas they
merge the scattered waves via the multiple scattering process among 
subdomains. (iv) We deal directly with solutions whereas they deal with
operator such as the scattering matrix or the Dirichlet to Neumann map.
\end{remark}
Consequently, we never solve a scattering problem in merging; they
solve an interesting and difficult scattering problem for each subdomain
and in every step of merging. 
\begin{definition}
\label{def-TW}
  A total wave in $D$ is a non-trivial solution of the homogeneous, 
  variable coefficient Helmholtz equation (\ref{2.1}) in $D$,
  subject to no boundary conditions. A total wave in a subdomain is a 
  non-trivial solution of (\ref{2.1}) in that subdomain, subject to no 
  boundary conditions, other than natural extension to the outside by 
  continuity conditions on the function and its normal derivative on
  the boundary.
\end{definition}
As such, it has a desirable, simple property: A total wave in 
a subdomain is utterly oblivious to the underlying multiple
scattering process among the subdomains. The total waves are much 
easier to construct and merge to build up the entire solution space for $D$.
\begin{observation}
\label{obs-merg2}
As two neighbor square subdomains merge, the two solution spaces, 
each of dimension $4L$, should merge to the solution space for the 
rectangle and of dimension $6L$. Indeed, the 2 continuity conditions 
(on function and its normal derivative) on the common interface of 
length $L$ of the two squares consist of $2L$ constraints on the 
$8L$ parameters of the original two solution spaces. 
\end{observation}
\begin{remark}
\label{rmk-dim1}
It is important to note that the neighbor subdomains may have different
$4L$ as dimensions of solution spaces. Similarly, the number $4L$ does
not imply that each edge of a square subdomain bears exactly 
$L$ parameters. In that case Observation \ref{obs-merg2} still holds:
Whatever number of parameters born on the common interface by the
Dirichlet, or Neumann, or D\&N data (see Observation \ref{obs-merg1})
will be eliminated for both subdomains after merging, and thus 
subtracted from the sum of dimensions of the two solution spaces. 
\end{remark}

\subsection{The algorithm}  \label{sec-3.2}

It should be noted that our new algorithm is almost indentical to
that of \cite{chen-FDS1} in structure and complexity; they differ 
in solution space and in the merging strategy. The complexities 
for the old and new differ by a constant multiple, with the new 
more efficient by a considerable factor. 

The main advantage of the new is its simplicity in discretization 
and merging, two major difficulties still being reckoned with in 
existing approaches. The total wave based new algorithm never solve 
any scattering problem in 
constructing the total wave solution spaces (TWSS) for subdomains and 
in merging the solution spaces. The new only deals with the total 
waves in subdomains until the solution space is obtained for the 
entire $D$, 
whereas the old methods solve a local scattering problem every step 
of the way, and have to deal with the typical corner singularities 
of the local scattered waves at the four corners of a subdomain. As 
two scattered wave solution spaces merge, four of the eight 
singularities at the eight corners of the two subdomains cancel out
and disappear. This cancellation would give rise to conditioning 
problems. The total waves in a subdomain have natural extensions 
to the outside. They are unaware of the manmade corners of the 
subdomains, and are easier to sample and merge.

We reiterate that for simplicity in describing our new algorithm, we 
assume that the scatterer $q$ is a smooth function in $D$, and 
vanishes smoothly outside square $D$. In actual implementation, we will 
relax this requirement to include piecewise smooth scatterers with jumps 
inside and on the boundary of a fairly arbitrary, bounded domain $D$ in 
two dimensions. We will also deal with jump discontinuities which form
corners in the original medium; see \S \ref{sec-3.3} for more details.
% f''(x) = [-1 16 -30 16 -1]/(12h^2) f(x) + O(h^4)
% f' (x) = [ 1 -8   0  8 -1]/(12h)   f(x) + O(h^4)

The total wave based fast direct solver is described in the following
four subsections.

\subsubsection{Construct TWSS for bottom level subdomains}  
\label{sec-3.2.1}
Step 1 of the algorithm is to construct total wave solution spaces
TWSS for each square subdomains $S$ at the bottom level of the quadtree. 

TWSS is obtained by solving the homogeneous, variable coefficient 
Helmholtz equation (\ref{2.1}) in $S$, subject to no boundary 
conditions. A high order method is provided in \S \ref{sec-3.3}, and
also implemented numerically. Here we present a second order method
to illustrate the procedure. The square $S$ is discretized $S$ a 
$m$-by-$m$ uniform mesh $Q$.
\begin{figure}[ht]
\begin{center}
\input{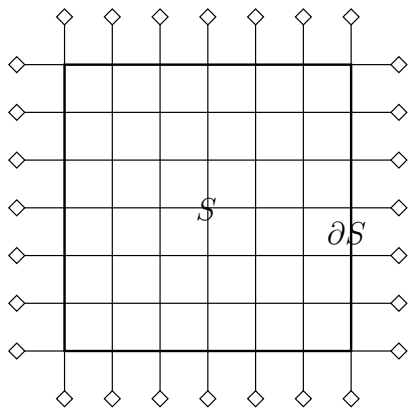}
\end{center}
  \vspace*{-0.1in}
  \caption{Uniform mesh on a square $S$}
  \label{fig-mesh2}
\end{figure}
Five point stencil, to replace the Laplacian of (\ref{2.1}) on every point, 
including those on the boundary $\partial S$, produces $m^2$ homogeneous 
equations
\begin{equation}
 \Delta_h u(x) + k^2 n(x) u(x) = 0,  \quad x \in Q
\label{3.8m}
\end{equation}
The equations for the $4(m-1)$ boundary points of the mesh will require
$4m$ additional, free variables $u(x)$, with $x$ a step $h$ away from the 
boundary; see Figure \ref{fig-mesh2} for points $x$ outside $S$ marked 
by $\Diamond$. All together, there are $m^2$ equations for $m^2 + 4m$ 
unknowns $u(x)$, and so the null space of the discrete operator 
$\Delta_h + k^2 n(x)$ is $4m$ dimensional, namely there
are $4m$ nontrivial solutions to the $m^2$ equations. 

These $4m$ basis functions for the TWSS in $S$ are collected in a matrix 
$U_s$ of size $m^2$-by-$4m$. The basis functions are also evaluate at the 
boundary $\partial S$, and when paired with their normal derivatives, 
provide the D\&N data for the TWSS. For simplicity, we assume that 
$\partial S$ is sampled with $4m$ points, and that there are $m$ points 
on each of the four edges of $S$, so that the matrices  
\begin{equation}
  U = U_s|_{\partial S} = \left[ u^{(1)}, u^{(2)}, \cdots, u^{(4m)} \right],  \quad
  U_n = \{\partial_n U_s\}|_{\partial S} = \left[ u_n^{(1)}, u_n^{(2)}, \cdots, u_n^{(4m)} \right ], \quad
\label{3.8a}
\end{equation}
are each of dimension $4m$-by-$4m$. Let 
\begin{equation}
  G = \left [ \begin{array}{l}
       U \\  U_n  
      \end{array} \right ], \quad
  \tilde{G} = \left [ \begin{array}{c}
       -U \\  U_n  
      \end{array} \right ]
\label{3.8b}
\end{equation}
so $G$ is of size $8m$-by-$4m$. Our numerical experiment shows that 
these $4m$ boundary points can be {\it equispaced}, as opposed to 
crowded toward the four corners. Denser sampling points are required 
near a singularity of the solution (arising from medium discontinuity), 
or near locations where the total waves have more evanescent modes due 
to medium complexity.
\subsubsection{Merging two subdomains}  
\label{sec-3.2.2}
Step 2 of the algorithm is bottom up merging: On each level of 
the quadtree, merge two neighbor square subdomains, whose TWSS's are 
available, to construct TWSS for the resulting rectangular domain. Merge 
again two neighbor rectangles to form TWSS for the resulting square domain 
on the higher level. Stop at the highest level which contains only one 
square domain that is $D$. Merging is achieved by imposing continuity 
conditions on the D\&N data on the common interface of two subdomains; 
see Observation \ref{obs-merg2}. 
\begin{figure}[ht]
\begin{center}
\input{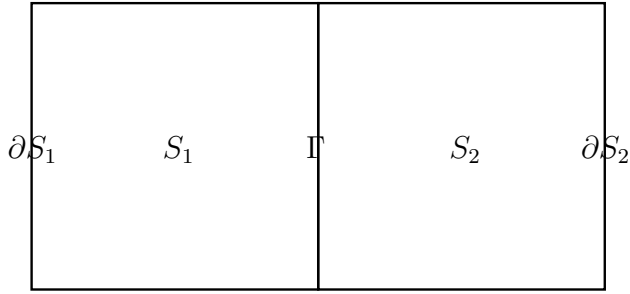}
\end{center}
  \vspace*{-0.1in}
  \caption{Merge two squares $S_1$ and $S_2$}
  \label{fig-mesh4}
\end{figure}
Let $G_1$ and $G_2$ be the D\&N data matrix of the two square subdomains $S_1$ 
and $S_2$. Let $\Gamma$ be their common interface; see Figure \ref{fig-mesh4}. 
Then merging $G_1$ and $G_2$ to produce D\&N data matrix $G$ for $S=S_1 
\cup S_2$ requires enforcing the continuity of the D\&N data across $\Gamma$, 
by solution of the homogeneous linear system of $2m$ equations
\begin{equation}
  G_1|_\Gamma \; \tau_1 = -\tilde{G}_2|_\Gamma \; \tau_2,  
\quad \mbox{namely} \quad 
  [ G_1|_\Gamma, \ \tilde{G}_2|_\Gamma ] \left[ \begin{array}{l}
                     \tau_1 \\  \tau_2  
                    \end{array} \right ] = 0
\label{3.9}
\end{equation}
where $G_i|_\Gamma$ is the D\&N data matrix $G_i$ restricted on $\Gamma$; 
it is a matrix of size $2m$-by-$4m$. The coefficient matrix $[G_1|_\Gamma, 
\ \tilde{G}_2|_\Gamma]$ is $2m$-by-$8m$ with a null 
space of dimension $6m$. Let $T_1$ and $T_2$ be matrices of size $4m$-by-$6m$ 
whose columns consist of the $6m$ solutions $\tau_1$ and $\tau_2$ of 
(\ref{3.9}). The D\&N data matrix $G$ on $\partial S = \{\partial S_1 
\setminus \Gamma\} \cup \{\partial S_2 \setminus \Gamma\}$ are given by
\begin{equation}
  G|_{\{\partial S_1 \setminus \Gamma\}} = G_1|_{\{\partial S_1 \setminus \Gamma\}} \; T_1,
  \quad \mbox{and} \quad 
  G|_{\{\partial S_2 \setminus \Gamma\}} = G_2|_{\{\partial S_2 \setminus \Gamma\}} \; T_2
\label{3.10}
\end{equation}
We refer to (\ref{3.9}), (\ref{3.10}) as the merging formulas, to be used
throughout the bottom-up merging process. 

In the remainder of this subsection we define and determine splitting. 
Let $u$ be a total wave - solution of the Helmholtz equation (\ref{2.1}) in 
$S=S_1 \cup S_2$. Let $g$, $g_1$, $g_2$ be the D\&N data for $u$ on 
$\partial S$, $\partial S_1$, $\partial S_2$, respectively. Thus, there 
exist coefficients $\gamma$, $\gamma_1$, $\gamma_2$ such that 
\begin{equation}
  g = G \gamma, \quad g_1 = G_1 \gamma_1, \quad g_2 = G_2 \gamma_2 
\label{3.10a}
\end{equation}
\begin{definition}
\label{def-split}
Splitting is the operation to determine $g_1, g_2$ from $g$, in terms of
their coefficients. The linear map 
$S_p: \gamma \mapsto
 \left[ \begin{array}{l}
     \gamma_1 \\ \gamma_2  
 \end{array} \right ] $ is referred to as the splitting operator.
\end{definition}
It follows (\ref{3.10}) immediately that $S_p$, of size $8m$-by-$6m$, is
given by the formula
\begin{equation}
  S_p = \left[ \begin{array}{l}
                     T_1 \\  T_2 
                    \end{array} \right ]
\label{3.10b}
\end{equation}
\begin{remark}
\label{rmk-bdry}
All the matrix operations after Step 1 are carried out on the boundaries
of the subdomains each with some $4m$ points on the boundary, instead of 
inside the subdomains each with about $m^2$ points. Therefore, the 
matrix operations are not nearly as costly. The entire merging step will
cost only $O(k^3)$ flops for a $k$-by-$k$ wavelength problem on $D$; see
\cite{chen-FDS1} for a complete analysis. 
\end{remark}
Subsections \ref{sec-3.2.1} and \ref{sec-3.2.2} are about constructing 
null spaces, and no scattering problem has been solved so far.
\subsubsection{Decompose TWSS via Green's formula}  
\label{sec-3.2.3}
Step 3 of the algorithm is to decompose TWSS for the whole
scatterer $D$ into the incoming and outgoing parts.
Now that the TWSS is constructed for $D$ with the D\&N data matrix $G$
available on $\partial D$ which lies in the free space, the Green's
identities apply, with the free space Green's function. In particular,
the D\&N data matrix $G$ for the total waves can be split by the projector 
$P_-$, see \S \ref{sec-4.2} for definition and technical details, to 
obtain the D\&N data matrix $G_0$ for the incident parts of the total waves. 

Let the incident components of $G$ be denoted by 
\begin{equation}
  G_0 = \left [ \begin{array}{r}
       U_0 \\  \partial_n U_0  
      \end{array} \right ] \quad \mbox{so that} \quad 
   G_0 = P_- \ G
\label{3.8d}
\end{equation}
Expressing the D\&N data of the prescribed incident wave $u_0$ on 
$\partial D$ by the basis $G_0$
\begin{equation}
  G_0 \ \gamma = \left [ \begin{array}{r}
       u_0 \\  \partial_n u_0  
      \end{array} \right ]   \label{3.8e}
\end{equation}
we solve (\ref{3.8e}) to obtain the coefficients $\gamma$. Obviously, 
\begin{equation}
  \left [ \begin{array}{r}
       u \\  \partial_n u  
      \end{array} \right ] =  G \ \gamma  \label{3.8f}
\end{equation}
is the D\&N data of the total wave corresponding to the prescribed 
incident wave $u_0$. For many applications, we want the scattered 
wave, and perhaps also its normal derivative, on $\partial D$. In this
case, the D\&N data of the scattered wave is obtained by subtracting
(\ref{3.8e}) from (\ref{3.8f}). This concludes our algorithm if the 
scattered wave on $\partial D$ is all we want; otherwise continue to
the next step. 
\subsubsection{Split a total wave}  
\label{sec-3.2.4}
Step 4 of the algorithm construct the total wave inside $D$. It is a 
top-down splitting process to propagate recursively along the quadtree 
the coefficient $\gamma$ from a domain $S$ to its subdomains $S_1$ and 
$S_2$. 

This is the reverse of Step 2 detailed in \S \ref{sec-3.2.2}, see 
(\ref{3.10b}) for details. Continue the recursive splitting till $\gamma$
is available for every bottom level square $S$. According to (\ref{3.8a}),
the total wave $u$ in $S$ is obtained by
\begin{equation}
  u = U_s \; \gamma \label{3.8n}
\end{equation}
Now the total wave $u$ of (\ref{3.8f}) is available everywhere inside $D$. 

This is the end of our algorithm, and the volume scattering problem 
for a given incident wave $u_0$ is solved.

\subsection{Some implementation details and remarks}  \label{sec-3.3}

There are several very accurate methods to construct TWSS for a small 
subdomain $S$ on the bottom level of the quadtree. We will present two
typical methods: Collocation and the weak formulation.

\subsubsection{Collocation method to construct TWSS}  
\label{sec-3.3.1}
The solutions $u$ of (\ref{2.1}) in $S$ can be approximated by polynomials
or bandlimited functions or some other suitable basis $B_j(x)$. Thus 
we represent the total waves $u$ in $S$ by
\begin{equation}
   u(x) = \sum_j c_j B_j(x), \quad x \in S
\label{3.11}
\end{equation}
The homogeneous, variable coefficient Helmholtz equation (\ref{2.1}), 
with $u$ given above, is evaluated at some $n$ suitable locations in 
$S$, giving rise to $n$ homogeneous linear equations for $c_j$
\begin{equation}
   A c = 0   \label{3.11a}
\end{equation}
The vectors $c$ from the null space of matrix $A$ are then used to 
construct the total waves $u$ for the TWSS in $S$. 

In numerical experiments, we used bandlimited functions of the form
\begin{equation}
  B(k,\theta,x) = \exp(ik(x_1\cos(\theta)+ x_2\sin(\theta))), \quad x \in S,
    \quad \theta \in [0, 2\pi), \quad k \in [k_1, k_2]   \label{3.11b}
\end{equation}
The $n$ collocation points on $S$ are either uniform mesh - see 
Figure \ref{fig-mesh2} - or a graded mesh such as the tensor 
legendre points.

\subsubsection{Weak formulation for TWSS}  
\label{sec-3.3.2}
The variational formulation for the Helmoltz equation (\ref{2.1}) offers
more flexibility for an irregular subdomain $S$, and leads to typical
finite element solution for (\ref{2.1}). Without boundary conditions
for (\ref{2.1}), the weak formulation assumes the form
\begin{equation}
-\int_S \nabla v \cdot \nabla u dx  +  \int_{\partial S} v u_n ds  + 
\int_S k^2(x) v(x) u(x) dx  = 0
 \label{3.11c}
\end{equation}
where $u$ and $v$ are both from a function space such as the one
spanned by basis in (\ref{3.11}), or by a typical finite element basis. 
In numerical experiments, we used (\ref{3.11b}) as basis.

\subsubsection{Sampling and completeness of TWSS}  
\label{sec-3.3.3}
For the prescribed precision $\epsilon$, the completeness of TWSS
depends only on sampling, specifically the number of points and
their distribution in each bottom level subdomain. Sampling rate
is determined, as is well known, by the local wave number, by the
complexity of the local medium, and by the distance to the nearest
singularities arising from (i) the sources of incident wave (ii) 
the corners (iii) the edges. These sources have different strength
of singularities, but if $\epsilon$ is small it is hardly necessary
to treat them differently in sampling rate near them. In the 
standard case when the medium changes smoothly and slowly relative 
to wavelength, and if the subdomain is far from a singularity, 2 
to 4 points per local wavelength is usually sufficient.

On the other hand, since we always over sample, the TWSS will be
unavoidably over complete to the prescribed precision. Once the 
TWSS space is constructed for a subdomain, particularly the bottom
level ones, it may be compressed by SVD, or by pivoted Gram-Schmidt 
or QR. This will only make the TWSS healthier (because the total 
waves will be orthonormalized). It will not reduce the level of
over-completeness (because smaller singular values of matrix $G$
of (\ref{3.8b}) may not correspond to more evanescent waves).

\subsubsection{Discontinuities in the medium}  
\label{sec-3.3.4}
If a bottom level subdomain $S$ contains a smooth interface across 
which the index of refraction jumps, then $S$ must be divided into 
two subdomains. These subdomans of $S$ becoms bottom level subdomains.

Let's assume now that a bottom level subdomain $S$ contains a corner 
of the medium over which the index of refraction jumps. Then $S$
must be partitioned into two subdomains, one contains the corner, the
other is the complement. TWSS for each subdomain will contain some 
regular solutions and some singular solutions because of the corner. 
The regular ones will be constructed, say, by the collocation method 
(\ref{3.11}). The singular ones will also be constructed by (\ref{3.11}), 
except that the basis functions for them must be chosen to contain the 
local singular behavior of the total wave near the corner. Finally, the 
regular and singular solutions are collected together on the boundary 
and compressed by SVD or QR.
 
As is well established, near such a corner the solution is spanned
by the Bessel-Fourier terms of fractional orders,
\begin{equation}
   u(x) \sim \sum_\nu c_\nu J_\nu(kr) e^{i\nu \theta}
\label{3.31}
\end{equation}
where $x$ near the corner is assigned a polar coordinates $(r,\theta)$
centered at the corner.

\subsubsection{Equations v.s. solutions}  
\label{sec-3.3.5}
The existing numerical methods for scattering or general elliptic problems 
can be divided into three categories, according to how much they are 
involved in building solution space. 
\begin{enumerate}
\item Our total wave approach deals with the total waves.
The differential equations are treated only at the bottom level
subdomains, and without boundary conditions. Merging TWSS of 
subdomains requires no knowledge of the PDEs.

\item Methods based on scattering matrix or related objects 
\cite{chen-FDS1}-\cite{h1-mat1} deal with the scattered waves for 
each subdomain on every level of the hierarchical domain 
decomposition (such as a quadtree). These methods solve a scattering 
problem on each bottom level subdomain, and they also need to know 
the underlying Green's function in merging two subdomains to 
construct the scattering matrix for the union of the two subdomains.
Each subdomain on the hierarchy has manmade corners, which is visible
to the scattered waves and show up as singularities in the 
scattering matrix. 

\item Finite difference, finite elements, or similar methods deals 
with the differential equations or their variational forms. Continuity
conditions over an interface of subdomains are enforced on the equations
rather than directly on the solutions, as our method does.
\end{enumerate}
Our total wave approach merges the subdomains directly and cleanly.
Existing methods have to avoid dealing directly with the unwieldy 
analytical issues around the manmade corners arising from the 
manmade subdomains, by merging two ``discretized'' subdomains, or
two overlapping ones.

\subsubsection{All merging formulas}  
\label{sec-3.3.6}

{\bf 1. Merging two squares $S_1, \ S_2$ to a rectangle $S$.} Let $G_1$ 
and $G_2$ be the two D\&N data matrices to be merged to produce the D\&N 
data matrix 
$G$ for $S=S_1 \cup S_2$. Denote by $G_{ij}$ the part of $G_j$ restricted 
to the edge shared with $G_i$, $i \not = j$. Then the merging-splitting 
matrices $T_1,\ T_2$ (see (\ref{3.9})) are solutions of the equation
\begin{equation}
   [ G_{21} \ \tilde{G}_{12}]_{2\times 8} \left[ \begin{array}{l}
                     T_1 \\  T_2  
                    \end{array} \right]_{8\times 6} = 0_{2\times 6}
\label{3.9r}
\end{equation}

{\bf 2. Merging two rectangles $S_1, \ S_2$ to a square $S$.} Let $G_1$ 
and $G_2$ be the two D\&N data matrices to be merged to produce the D\&N 
data matrix 
$G$ for $S=S_1 \cup S_2$. Denote by $G_{ij}$ the part of $G_j$ restricted 
to the edge shared with $G_i$, $i \not = j$. Then the merging-splitting 
matrices $T_1,\ T_2$ are solutions of the equation
\begin{equation}
   [ G_{21} \ \tilde{G}_{12}]_{4\times 12} \left[ \begin{array}{l}
                     T_1 \\  T_2  
                    \end{array} \right]_{12\times 8} = 0_{4\times 8}
\label{3.9s}
\end{equation}

{\bf 3. Merging four squares $S_i,\ i=1:4$ to a square $S$.} Let $G_i,\ 
i=1:4$ be the four D\&N data matrices to be merged to produce the D\&N 
data matrix $G$ for $S=\cup_{i=1}^4 S_i$. Denote by $G_{ij}$ the part of 
$G_j$ restricted to the edge shared with $G_i$, $i \not = j$. Then the 
merging-splitting matrices $T_i,\ i=1:4$ are solutions of the equation
\begin{equation}
 \left[ \begin{array}{cccc}
   G_{21} & \tilde{G}_{12}  & 0              & 0 \\
   0     & G_{32}          & \tilde{G}_{23}  & 0 \\
   0     & 0              & G_{43}          & \tilde{G}_{34}  \\
   \tilde{G}_{41} & 0      & 0              & G_{14}  
 \end{array} \right]_{8\times 16} \left[ \begin{array}{c}
                     T_1 \\ T_2  \\ T_3  \\ T_4  
                    \end{array} \right]_{16\times 8} = 0_{8\times 8}
\label{3.9t}
\end{equation}

\subsection{Extensions to surface scattering problems}  \label{sec-3.4}
Our total wave approach for VSP can be extended to a surface scattering 
problem (SSP). For simplicity, we consider a 2-D SSP off a sound soft 
(zero Dirichlet for total wave) smooth convex body $D$ such as the
unit disc with  the sources of the incident wave $u_0$ well separated 
from $D$. Given $u_0$ on $\partial D$, the SSP is to determine 
$\partial_n v$, the normal derivative of the scattered wave $v$ on 
$\partial D$.

\noindent The following steps outline a possible extension of the TWSS
method to SSP.
\begin{enumerate}
\item {\it The bottom level subdomains.}
Dividing the annulus $1\le r \le 1+h$, for some $h\ge 0$, along
the radial direction into sufficient many pieces. Each piece 
\begin{equation}
   A_i = \{(r,\theta), \; 1\le r \le 1+h, \; \theta_{i-1} \le
             \theta \le \theta_i \}, \quad i=1:n
\label{3.36}
\end{equation}
is bounded by four curves: two arcs and two straight radial segments. 
Remove the arc of $r=1+h$. The remaining three curves form the $i$-th 
subdomain $\Gamma_i$ on the bottom level:
\begin{eqnarray}
   \Gamma_{i1} &=& \{(1,\theta), \; \theta_{i-1} \le \theta \le 
                   \theta_i \} \label{3.37a}  \\
   \Gamma_{i2} &=& \{(r,\theta_{i-1}), \; 1\le r \le 1+h \} \label{3.37b}  \\
   \Gamma_{i3} &=& \{(r,\theta_i), \; 1\le r \le 1+h \} \label{3.37c}   \\
   \Gamma_i &=& \Gamma_{i1} \cup \Gamma_{i2} \cup \Gamma_{i3} \label{3.37}
\end{eqnarray}
\item {\it TWSS in bottom level subdomains.} For a prescribed
  precision construct TWSS in $A_i$ subject to the zero Dirichlet
  condition (sound soft) on the arc $\Gamma_{i1}$. Let there
  be $m$ total waves in TWSS of the form $u^{(j)} =u^{(j)}_0 +
  v^{(j)}$, $j=1:m$. For each $u^{(j)}$, the TWSS will contain four 
  functions (i) $u^{(j)}$ and $\partial_n u^{(j)}$ restricted on the 
  two radial segments $\Gamma_{i2}, \Gamma_{i3}$ (ii) $u^{(j)}_0$ and 
  $\partial_n v^{(j)}$ on the arc $\Gamma_{i1}$.
\item {\it Merging.} Merge the subdomains, starting from the bottom
  level ones $\Gamma_i$, recursively and upward along a hierarchical 
  domain decomposition. Two neighboring subdomains $C_1, \, C_2$, 
  separated by a common radial segment, will be merged to form their 
  parent submain $P$ by continuities of total wave and its normal 
  derivative on the interface.
\item {\it After merging.} The TWSS will contain pairs of functions
on $\partial D$, the unit circle. Each pair is of the form
\begin{equation}
   u^{(\ell)}_0, \; \partial_n v^{(\ell)}, \quad \ell=1:N
\label{3.38}
\end{equation}
where $\partial_n v^{(\ell)}$ is the normal derivative of the
scattered wave $v^{(\ell)}$ off $D$ induced by an incident wave 
$u^{(\ell)}_0$; $N$ is proportional to arclength of $\partial D$ 
measured in wavelength. 

\item {\it Construct $\partial_n v$.} Finally we construct the 
normal derivative, on $\partial D$, of the scattered wave $v$ 
off $D$ induced by the prescribed incident wave $u_0$. Spanning 
$u_0$ by $u^{(\ell)}_0$
\begin{equation}
   u_0 = \sum_\ell \alpha_\ell u^{(\ell)}_0
\label{3.39}
\end{equation}
we use the coefficients to produce $\partial_n v$
\begin{equation}
   \partial_n v = \sum_\ell \alpha_\ell \partial_n v^{(\ell)}
\label{3.40}
\end{equation}
\end{enumerate}
This is the end of the algorithm. We have constructed the pair
$(v,\partial_n v)$ on $\partial D$, with $v=-u_0$, and solved
the SSP for sound soft scatterer $D$. The solution of (\ref{3.39}) 
for $\alpha_\ell$ is a bottleneck of the procedure.

\section{The Volume Scattering Problem}  \label{sec-4}

One of the misfortunes of the 20th Century applied mathematics
is that the volume scattering problem (VSP) for the Helmholtz 
equation was posed and is still being solved today as a boundary 
value problem for the scattered wave in a domain $D$, and worse
yet, for each subdomain of $D$ in a domain decomposition setting. 
The scattering problem is a very special boundary value problem 
in that both the Dirichlet and Neumann data of the incident wave 
are available on the boundary of the whole scatterer $D$. In this 
section we present the classical formulations for VSP, and special 
properties useful for the total wave based fast direct solver. 

Given index of refraction $n(x) = \sqrt{1+q(x)}$ in a bounded domain $D$,
we consider the volume scattering problem in $k$-space governed by the 
Helmholtz equation
\begin{equation}
  \Delta u + k^2 n^2(x)u = 0,   \quad \mbox{or}  \quad
  \Delta u + k^2(1+q)  u = 0 
\label{2.1}
\end{equation}
where $u$ is the total wave, $q$ is the scatterer, $n=1$ and $q=0$ 
outside $D$. In a typical setting, $u$ is of the form 
\begin{equation}
  u = u_0 + v,
\label{2.2}
\end{equation}
where the incident wave $u_0$ is given in $D$ and the scattered wave 
$v : \R^2\mapsto \C$ is to be determined as a solution of the 
{\it inhomogeneous} Helmholtz equation, 
\begin{equation}
  \Delta v + k^2 v = -k^2 q(u_0+v)
     \label{2.3}
\end{equation}
subject to the Sommerfeld radiation condition
\begin{equation}
   \lim_{r\rightarrow \infty} \sqrt{r}
        \left( \frac{\partial v}{\partial r} -
        ik v \right) = 0  \label{2.3a}
\end{equation}
The scattering problem (\ref{2.3}), (\ref{2.3a}) can also be
formulated as the Lippmann-Schwinger equation
\begin{equation}
   \sigma(x) + k^2 q(x) \int_D G(x,\xi) \sigma(\xi) d\xi  
 = -k^2 q(x) u_0(x) \label{3.17}
\end{equation}
for the monopole density $\sigma$ in $D$, which is related to
$v$ by
\begin{equation}
   v(x) =  \int_D G(x,\xi) \sigma(\xi) d\xi.  \label{3.18}
\end{equation}
where $G = -(i/4) H_0(k|x-\xi|)$ is the free space Green's function.

\subsection{Dirichlet and Neumann Data for $u_0$ on $\partial D$}  
\label{sec-4.1}

For the scattering problem (\ref{2.3}) and (\ref{2.3a}), or 
(\ref{3.17}) and (\ref{3.18}), the incident wave $u_0$ must be
available inside the scatterer $D$, or better, its sources outside
$D$ are prescribed.
\begin{observation}
\label{obs-DNdata}
For a bounded domain $D$ with a regular boundary $\partial D$, the 
Dirichlet and Neumann data for the incident wave $u_0$ are always
available on $\partial D$.
\end{observation}
To verify this statement, we observe that more often then not in a 
typical application, the incident wave $u_0$ of a volume scattering 
problem is specified by its sources outside $D$, such as a monopole
or a plane wave. In that case, the D\&N data $u_0$ and $\partial_n
u_0$ can be evaluated directly on the boundary. 

Suppose $u_0$ is only given in $D$, as required by (\ref{2.3}) or
(\ref{3.17}). By Green's third identity, $u_0$, being solution of 
$\Delta v + k^2 v = 0$ in $D$, is given by its D\&N data
\begin{equation}
 u_0(x) = -\int_{\partial D} \left( \partial_n u_0(\xi) \cdot 
          G(x,\xi) - u_0(\xi) \cdot \frac{\partial G(x,\xi)}{\partial 
                n(\xi)} \right) ds(\xi),      \label{3.13}
\end{equation}
So the D\&N data for $u_0$ can be recovered by solving (\ref{3.13}) 
as an equation. This solution process can often be simplified if
$u_0$ is given not only in $D$, but also on $\bar{D}$.

\subsection{Green's Third Identity as the Interior Projector $P_-$}  
\label{sec-4.2}

Let $W = L^2(\partial D) \times L^2(\partial D)$. Let $W_\pm$ be 
two subspaces of $W$ defined by the formulae
\begin{eqnarray}
   W_- &=& \left\{ \ \left. \left( v, \partial_n v 
           \right) \ \right| \ \Delta v + k^2 v = 0 \ \mbox{in} \ D 
           \ \right\}      \label{1.2a}  \\
   W_+ &=& \left\{ \ \left. \left( v, \partial_n v 
           \right) \ \right| \ \Delta v + k^2 v = 0 \ 
                   \mbox{outside} \ \bar{D} \
           \mbox{subject to (\ref{2.3a})} 
           \ \right\}      \label{1.2b}
\end{eqnarray}
Therefore, $W_-$ consists of the D\&N data for incident waves
in $D$, and $W_+$ consists of the D\&N data for outgoing (or 
scattered) waves outside $\bar{D}$. 

For a bounded domain $D$ with a regular boundary $\partial D$, the 
linear map $(\phi, \psi) \in W \mapsto v \in L^2(D)$ defined by
Green's third identity
\begin{equation}
 v(x) = -\int_{\partial D} \left( \partial_n \phi(\xi) \cdot G(x,\xi) 
                   - \phi(\xi) \cdot \frac{\partial G(x,\xi)}{\partial 
                n(\xi)} \right) ds(\xi), 
        \label{3.13a}
\end{equation}
produces $v$ with $\Delta v + k^2 v = 0$ in $D$. Let $x$ approach
$\partial D$ from $D$, and denote the limit by $x_-$ corresponding 
to $x \in \partial D$. We thus obtain the linear map $P_-: (\phi,
\psi) \in W \mapsto (v, \partial_n v) \in W_-$
\begin{eqnarray}
 v(x) &=& -\int_{\partial D} \left( \partial_n \phi(\xi) \cdot G(x_-,\xi) 
                   - \phi(\xi) \cdot \frac{\partial G(x_-,\xi)}{\partial 
                n(\xi)} \right) ds(\xi),        \label{3.13b} \\
 \partial_n v(x) 
      &=& -\int_{\partial D} \left( \partial_n \phi(\xi) 
                \cdot \frac{\partial G(x_-,\xi)}{\partial n(x)} 
           - \phi(\xi) \cdot \frac{\partial^2 G(x_-,\xi)}{\partial n(x)
             \partial n(\xi)} \right) ds(\xi),        \label{3.13c} 
\end{eqnarray}
Likewise, using Green's third identity outside $D$, we introduce
another linear map $P_+: (\phi,\psi) \in W \mapsto (v, \partial_n v) 
\in W_+$. In terms of the standard layer potential operators $S, K, K', T$ 
defined by (\ref{7.26})-(\ref{7.29}), the two operators $P_\pm$ are 
given by
\begin{eqnarray}
   P_+
   \left [ 
     \begin{array}{c} 
       \phi \\
       \psi
     \end{array}
   \right] 
=  \left\{ 
    \frac{1}{2}
    \left[ 
     \begin{array}{cc} 
       I  &  0 \\
       0  &  I
     \end{array}  
    \right]
+ 
    \left[ 
     \begin{array}{cc} 
       -K  &  S \\
       -T  &  K'
     \end{array}  
    \right] 
   \right\}
   \left [ 
     \begin{array}{c} 
       \phi \\
       \psi
     \end{array}
   \right] 
      \label{3.14p}\\
   P_-
   \left [ 
     \begin{array}{c} 
       \phi \\
       \psi
     \end{array}
   \right] 
=  \left\{ 
    \frac{1}{2}
    \left[ 
     \begin{array}{cc} 
       I  &  0 \\
       0  &  I
     \end{array}  
    \right]
- 
    \left[ 
     \begin{array}{cc} 
       -K  &  S \\
       -T  &  K'
     \end{array}  
    \right] 
   \right\}
   \left [ 
     \begin{array}{c} 
       \phi \\
       \psi
     \end{array}
   \right] 
   \label{3.14q}
\end{eqnarray}
It follows immediately from Green's third identity (\ref{3.13}) that
\begin{observation}
\label{obs-proj}
The linear map $P_\pm : W \mapsto W_\pm$ is a projector converting
an arbitrary pair of boundary data $(\phi, \psi) \in W$ to the D\&N
data for an outgoing/incident wave. In particular, $P_-$ maps the D\&N 
data of the total wave $u$ of (\ref{2.2}) to those of its incident 
component $u_0$, and $P_+$ maps the D\&N data of total wave $u$ to 
those of its scattered component $v$
\end{observation}
\begin{remark}
\label{rmk-Pfirst}
Step 3 of the algorithm splits the D\&N data of a total wave on 
$\partial D$ for its scattered and incident parts using $P_-$. 
Since the D\&N data is nearly twice redundant, we 
only need the D\&N data for the incident part evaluated at about 
half as many discretization points (quadrature nodes for $P_-$).
We will select these points away from corner and edges to avoid
singularities of $P_-$ there.
\end{remark}
In addition, if the wave number $k$ is not a Dirichlet eigenvalue of 
the Laplacian in $D$, namely if $k$ does not hit a geometric resonance
of $D$, then only the first half of $G_0$ in (\ref{3.8e}) is required
to determine $\gamma$. In other words, the second half of $P_-$, which
involves the hyper singular kernel $T$, is not necessary for splitting.

\section{Layer potential representation for $P_+, P_-$}
\label{sec-5}

For $x$ close to the the boundary $\partial D$, the single and double
layer potentials
\begin{eqnarray}
   p(x) &=& \int_{\partial D} G(x,\xi) \psi(\xi) ds(\xi),  \label{7.22a} \\
   q(x) &=& \int_{\partial D} \frac{ \partial G(x,\xi) }{ \partial n(\xi) }  
             \phi(\xi)  ds(\xi)    \label{7.23a} 
\end{eqnarray}
associated with arbitrary pair of functions $(\phi,\psi)$, can be rewritten
\begin{eqnarray}
   p(x \pm hn(x)) &=& \int_{\partial D} G(x \pm hn(x),\xi)  
     \psi(\xi) ds(\xi),  \label{7.25a} \\
   q(x \pm hn(x)) &=& \int_{\partial D} \frac{ \partial G(x \pm hn(x),\xi) 
     }{ \partial n(\xi) } \phi(\xi) ds(\xi)    \label{7.26a} 
\end{eqnarray}
where $x$ now is on the boundary $\partial D$. Taking the directional
derivative of $p, q$ in the normal direction $n(x)$, we have
\begin{eqnarray}
    \frac{ \partial p(x \pm hn(x)) }{ \partial n(x) }
&=& \int_{\partial D} \frac{ \partial G(x \pm hn(x),\xi) }{ \partial n(x) }
   \psi(\xi) ds(\xi), \label{7.25e} \\
    \frac{ \partial q(x \pm hn(x)) }{ \partial n(x) } 
&=& \int_{\partial D} \frac{ \partial^2 G(x \pm hn(x),\xi) 
     }{ \partial n(x) \partial n(\xi) } \phi(\xi) ds(\xi)   
       \label{7.26e} 
\end{eqnarray}
Therefore, the Dirichlet and Neumann data $(\phi, \phi_n)$
of (\ref{3.13a}) are given by
\begin{equation}
   \left [ 
     \begin{array}{c} 
       \phi \\
       \phi_n
     \end{array}
   \right ]  
=
P_-
   \left [ 
     \begin{array}{c} 
       \phi \\
       \psi
     \end{array}
   \right]
=
   -\lim_{h\rightarrow +0}
   \left[ 
     \begin{array}{c} 
       p(x^-)-q(x^-) \\
       \frac{ \partial }{ \partial n(x) } (p(x^-) - q(x^-))
     \end{array}  \right ]
      \label{3.13i}
\end{equation}
for $x^-=x-hn(x)$, $x \in \partial D$. The Dirichlet and Neumann 
data $(\phi, \phi_n)$ of (\ref{3.13c}) are given by
\begin{equation}
   \left [ 
     \begin{array}{c} 
       \phi \\
       \phi_n
     \end{array}
   \right ]  
=
P_+
   \left [ 
     \begin{array}{c} 
       \phi \\
       \psi
     \end{array}
   \right] 
=
   \lim_{h\rightarrow +0}
   \left[ 
     \begin{array}{c} 
       p(x^+)-q(x^+) \\
       \frac{ \partial }{ \partial n(x) } (p(x^+) - q(x^+))
     \end{array}  \right ]
      \label{3.13e}
\end{equation}
for $x^+=x+hn(x)$, $x \in \partial D$. The use of the jump conditions
\begin{eqnarray}
   \lim_{h\rightarrow +0} q(x \pm hn(x)) 
&=& \int_{\partial D} \frac{ \partial G(x,\xi) }{ \partial n(\xi) }  
     \phi(\xi) ds(\xi) \mp \frac{1}{2} \phi(x),   \label{7.24} \\
   \lim_{h\rightarrow +0} \frac{ \partial p(x \pm hn(x)) }{ \partial n(x) }
&=& \int_{\partial D} \frac{ \partial G(x,\xi) }{ \partial n(x) }  
        \psi(\xi) ds(\xi) \pm \frac{1}{2} \psi(x),   
                \label{7.25} 
\end{eqnarray}
reduces (\ref{3.13i}) and (\ref{3.13e}) to expressions for $P_\pm$
\begin{equation}
   P_\pm = \frac{1}{2}
    \left[ 
     \begin{array}{cc} 
       I  &  0 \\
       0  &  I
     \end{array}  
    \right]
\pm 
    \left[ 
     \begin{array}{cc} 
       -K  &  S \\
       -T  &  K'
     \end{array}  
    \right] 
   \label{3.14i}
\end{equation}
where the layer potential operaters $S, K, K', T$ are defined by
\begin{eqnarray}
   (S\sigma)(x) &=& \int_{\partial D} G(x,\xi) \cdot \sigma(\xi) \cdot ds(\xi),  
                                                 \label{7.26} \\
   (K\sigma)(x) &=& \int_{\partial D} \frac{ \partial G(x,\xi) }{ \partial 
     n(\xi) } \cdot \sigma(\xi) \cdot ds(\xi),   \label{7.27} \\
   (K'\sigma)(x) &=& \int_{\partial D} \frac{ \partial G(x,\xi) }{ \partial 
     n(x) } \cdot \sigma(\xi) \cdot ds(\xi),   \label{7.28} \\
   (T\sigma)(x) &=& \lim_{h\rightarrow +0} \int_{\partial D} 
     \frac{ \partial^2 G(x \pm hn(x),\xi) }{ \partial n(x) \partial 
     n(\xi) } \cdot \sigma(\xi) \cdot ds(\xi)   \label{7.29}
\end{eqnarray}
for $x \in \partial D$. Obviously, $P_+ + P_- = I$ implying that the
decomposition, of an arbitrary pair of data $(\phi,\psi)$ on boundary into
incoming and outgoing parts, is complete. 
\begin{remark}
\label{rmk-sigma}
It is well-known that for a smooth $\partial D$ the operators 
$S, \ K, \ K'$ are bounded from $C(\partial D)$ to $C(\partial D)$, 
whereas $T$ are bounded from $C^1(\partial D)$ to $C(\partial D)$. 
\end{remark}
Define the four linear operators $S_{\pm}, \ K_{\pm}, \ K'_{\pm}, \ T_{\pm}$
by the formulae
\begin{eqnarray}
  (S_{\pm}\psi)(x)  &=& \lim_{h\rightarrow +0} p(x \pm hn(x)), \label{7.23c} \\ 
  (K_{\pm}\phi)(x)  &=& \lim_{h\rightarrow +0} q(x \pm hn(x)), \label{7.24c} \\
  (K'_{\pm}\psi)(x) &=& \lim_{h\rightarrow +0} \frac{ \partial p(x \pm hn(x)) 
               }{ \partial n(x) },                       \label{7.25c} \\ 
  (T_{\pm}\phi)(x)  &=& \lim_{h\rightarrow +0} \frac{ \partial q(x \pm hn(x)) 
               }{ \partial n(x) }                        \label{7.26c}
\end{eqnarray}
for $x \in \partial D$, and we see that 
\begin{equation}
  K_{\pm}=K \mp \frac{1}{2}I, \quad  K'_{\pm}=K' \pm \frac{1}{2}I, \quad
  S_{\pm}=S, \quad  T_{\pm}=T.  
\label{7.15c}
\end{equation}
As is well-known, the exterior Green's formula maps to zero the Dirichlet 
and Neumann data $(\phi, \phi_n)$ of a outgoing wave $\phi$ in $D$.
Conversely, the interior Green's formula maps to zero the Dirichlet 
and Neumann data $(\psi, \psi_n)$ of a outgoing wave $\psi$ outside
$D$; therefore,  
\begin{equation}
   P_+ \circ P_- = P_- \circ P_+ =0, \quad W_- \cap W_+ = \{0\}.  
\label{7.15}
\end{equation}
\begin{lemma}
\label{lem-progen}
$P$ is a projector if and only if there exists a unique operator $Q$
such that 
\begin{equation}
  Q^2 = \frac{1}{4} I, \quad P = \frac{1}{2} I + Q
     \label{7.45}
\end{equation}
In particular,
\begin{equation}
   P_\pm = \frac{1}{2} I \pm Q, \quad \mbox{with} \quad
   Q = 
    \left[ 
     \begin{array}{cc} 
       -K  &  S \\
       -T  &  K'
     \end{array}  
    \right]
\end{equation}
Furthermore, $Q^2 = I/4$ implies that 
\begin{equation}
   S T = -\frac{1}{4} I + K^2, \quad  T S = -\frac{1}{4} I + (K')^2, \quad 
   K S = S K', \quad  T K = K' T
\end{equation}
Finally, the Dirichlet-to-Neumann maps $\Lambda_\pm: \phi \mapsto 
\phi_n$, with $(\phi, \phi_n) \in W_\pm(\partial D)$, are given by 
the formulae
\begin{eqnarray}
   \Lambda_\pm &=& \left(K' \mp \frac{1}{2} I \right)^{-1} T 
               = S^{-1} \left(K \pm \frac{1}{2} I \right)  \\
               &=& T \left(K \mp \frac{1}{2} I \right)^{-1} 
               = \left(K' \pm \frac{1}{2} I \right) S^{-1} 
\end{eqnarray}
\end{lemma}

\end{document}